\documentclass[11pt,a4paper]{article}
\usepackage{lipsum,url}
\usepackage[ruled]{algorithm2e}
\usepackage{float}
\usepackage{amssymb}
\usepackage{sidecap}
\usepackage{wrapfig}
\usepackage{tikz}
\usepackage{array}
\usepackage{mathptmx}
\usepackage{hyperref,xcolor}
\usepackage{caption}
\usepackage[T1]{fontenc}
\usepackage{libertine}
\usepackage{eurosym}
\usepackage{enumerate}
\usepackage{amsmath}
\usepackage{amssymb}
\usepackage{amsfonts}
\usepackage{amsthm}
\usepackage{amsfonts}
\usepackage{enumerate}
\usepackage{microtype}
\usepackage[english]{babel} 
\usepackage[utf8]{inputenc}
\usepackage{caption}
\usepackage{graphicx}
\usepackage{mathrsfs}
\usepackage{calrsfs}
\usepackage{pgf, tikz}
\usepackage{tikz-cd}
\usepackage{hypcap}
\usetikzlibrary{arrows,automata}
\graphicspath{{pictures/}} 
\usepackage{youngtab} 
\usepackage{ytableau} 
\expandafter\let\expandafter\oldproof\csname\string\proof\endcsname
\let\oldendproof\endproof
\renewenvironment{proof}[1][\proofname]{%
	\oldproof[\ttfamily \scshape #1]%
}{\oldendproof}
\usepackage[all]{xy}
\usepackage{mathrsfs}
\usepackage{tikz}

\usetikzlibrary{shapes,arrows}
\usetikzlibrary{plotmarks}
\usepackage{pgfplots}
\usepackage{wasysym}
\usepackage{tkz-graph}
\usepackage{tkz-berge}
\usetikzlibrary{shapes}
\usetikzlibrary{patterns}
\usetikzlibrary{matrix, arrows,topaths}
\usetikzlibrary{decorations.pathreplacing,decorations.markings}
\usetikzlibrary{decorations.markings,arrows}
\usetikzlibrary{decorations.pathreplacing,decorations.markings}
\tikzset{
	on each segment/.style={
		decorate,
		decoration={
			show path construction,
			moveto code={},
			lineto code={
				\path [#1]
				(\tikzinputsegmentfirst) -- (\tikzinputsegmentlast);
			},
			curveto code={
				\path [#1] (\tikzinputsegmentfirst)
				.. controls
				(\tikzinputsegmentsupporta) and (\tikzinputsegmentsupportb)
				..
				(\tikzinputsegmentlast);
			},
			closepath code={
				\path [#1]
				(\tikzinputsegmentfirst) -- (\tikzinputsegmentlast);
			},
		},
	},
	mid arrow/.style={postaction={decorate,decoration={
				markings,
				mark=at position .5 with {\arrow[#1]{stealth}}
			}}},
		}
		\usepackage{enumitem}
		\usepackage{enumerate}
		\theoremstyle{definition}
		\newtheorem{theorem}{Theorem}[section]
		\newtheorem{lemma}[theorem]{Lemma}
		\newtheorem{proposition}[theorem]{Proposition}
		
		\newtheorem{remark}[theorem]{Remark}
		\newtheorem{definition}[theorem]{Definition}
		\newtheorem{example}[theorem]{Example}

		\makeatletter
		\def\@part[#1]#2{%
			\ifnum \c@secnumdepth >-2\relax
			\refstepcounter{part}%
			\addcontentsline{toc}{part}{\thepart\hspace{1em}#1}%
			\else
			\addcontentsline{toc}{part}{#1}%
			\fi
			\markboth{}{}%
			\reset@font
			\parindent \z@ 
			\vspace*{10\p@}%
			\hbox{%
				\vbox{%
					\hsize=7mm%
					\begin{tabular}{@{}p{7mm}@{}}
						\makebox[7mm]{\scshape\strut\small\partname}\\
						\makebox[7mm]{\cellcolor{black}\Huge\color{white}\bfseries\strut\thepart\rule[-4cm]{0pt}{4cm}}%
					\end{tabular}%
					\makebox(0,0){\put(-10,-100){\fbox{\phantom{\rule[-4cm]{7mm}{4cm}}}}}
				}%
				\kern-2pt
				\vbox to 0pt{%
					\tabular[t]{@{}p{1cm}p{\dimexpr\hsize-2.1cm}@{}}\hline
					& \Huge\itshape\rule{0pt}{1.5\ht\strutbox}#1\endtabular}%
			}%
			\cleardoublepage
		}
		\makeatother
		\makeatletter
		\newcommand\HUGE{\@setfontsize\Huge{58}{67}} 
		\makeatother
		
		\usepackage{tabulary}
		
		\graphicspath{{graphics/}} 

\title{Two-Sided Noncommutative Gr\"{o}bner Basis on Quiver Algebras}
\author{Daniel K.  Waweru\\
	Ph.D  Student\\
		Nanjing University of Aeronautics and Astronautics\\
		Department of Communication and Information Engineering\\
		College of Electronic and Information Engineering\\
		29 Nanjing, Jiangsu 210016, China\\
		Email: daniel.kariuki19@yahoo.com\\	
	\and 
	Damian M Maingi, Ph.D. \\
			Assistant Professor\\
			Sultan Qaboos University\\
			Department of Mathematics, RM 0228\\
			College of Science, PB 36\\
			Muscat 123, Oman\\
			Email: dmaingi@squ.edu.om
	}	
\begin{document}
	\maketitle
	\begin{abstract}
	For a quiver $Q$, we define a path algebra $KQ$ as a span of all the paths of positive length. We study left (respective right) sided ideals and their Gr\"{o}bner bases. We introduce the two-sided ideals, two-sided division algorithm for elements of $KQ$ and study the two-sided Gr\"{o}bner bases. We show that with the defined two-sided division algorithm and two-sided Buchberger's algorithm, we can find a finite or an infinite Gr\"{o}bner basis for a two-sided ideal $I \subseteq KQ$ given a fixed admissible ordering. 
	\end{abstract}
	
	\section{Introduction}
	
	In 1986, Teo Mora published a paper \cite{mora} giving an algorithm for constructing a noncommutative Gr\"{o}bner Basis. This work built upon the work of George Bergman in particular his diamond lemma for ring theory \cite{bergman}. Mora's algorithm and the theory behind it, in many ways gives a noncommutative version of the Gr\"{o}bner Basis theory as seen in the commutative case which states; given an initial set $F$ generating an ideal $I$ in a polynomial ring $A$, Gr\"{o}bner basis theory uses $F$ to find a basis $G$ for $I$ with the property that for any $f \in  A$, division of $f$ by $G$ has a unique remainder. \\
	
	How we obtain that Gr\"{o}bner Basis remains the same as in commutative case where we add nonzero S-polynomials to an initial basis. The difference comes in the definition of an S-polynomial.  Since the purpose of S-polynomial $S(f,g)$ for each pair of nonzero polynomials $f,g \in A$ is to ensure that any polynomial $h \in A$ reducible by both $f$ and $g$ has a unique remainder when divided by a set of polynomials containing both $f$ and $g$, in commutative case there is only one way to divide $h$ by $f$ and $g$ giving the reduction $(h-x_1f)$ or $(h-x_2g)$ respectively, where $x_1$ and $x_2$ are terms. Thus there is only one S-polynomial for each pair of polynomials.	In noncommutative case however, a polynomial may divide another in many different ways. Therefore we do not have a fixed number of S-polynomials for each pair of polynomials in $A$. The number of S-polynomials depend on the number of overlaps between the leading monomials of $f$ and $g$. We can therefore strengthen division algorithm axiomatically to eliminate much ambiguity.
	
	\section{Preliminary}
We start off by mentioning some of rudiments to noncommutative Gr\"{o}bner basis over a noncommutative polynomial ring. 
	\subsection{Noncommutative Gr\"{o}bner Basis in Polynomial Ring}\label{sec2.1}
In this subsection $A=K[x_1, \dots, x_n]$ is a noncommutative polynomial ring. Monomials of $A$ are generated by alphabetical words of $A$ over K. We denote set of all monomials of $A$ by $M$. We thus make the following important definitions.
		\begin{definition}\label{def2.2} A relation $\prec$ is said to be a noncommutative monomial ordering on set $M$ if it satisfies;
		\begin{enumerate}[noitemsep]
			\item $\prec$ is a total order on $M$.
			\item $x\succ 1, \forall x \in M$.
			\item $x \succ y \Rightarrow wxz \succ wyz, \quad \forall x,y,w,z \in M.$
		\end{enumerate}
	\end{definition}
\begin{definition} Let $x,y \in M$. An $x-y$ overlap occurs when one can find factors $x=x_1z, y=zy_1$ where $x\neq x_1$ and $y\neq y_1.$  Different factorization in $M$ gives different overlaps.
\end{definition}
	\begin{definition}\label{def2.5} Every element $f \in A$ has a unique form $f=\sum_{i=1}^{n} a_iw_i, \quad a_i \in K, \quad w_i \in \langle x_1, \dots, x_n\rangle.$ We will denote all such monomials $w_i$ appearing in $f$ by $Mon(f)$. Furthermore for a noncommutative monomial ordering $\prec$, $w$ is called the leading monomial of $f \in A$, denoted by $w=LM(f)$ if $w$ occurs in $f$ and $w \succ m$ for all monomials $m \in Mon(f).$ The coefficient of $LM(f)$ in $f$ is called the leading coefficient and is denoted by $LC(f).$ The term $LT(f)=LC(f)LM(f)$ is called the leading term of $f$. Let $J \subset A$ then we define $LT(J)=\lbrace LT(g) :  g \in J \rbrace$.
	\end{definition}
		 For nonzero polynomials $f,g \in A$, we say that $f$ divides $g$ if the leading term of $f$ divides some term $h$ in $g$, where $h=x_l LM(f)x_r$ for some monomials $x_l$ and $x_r$. For noncommutative cases, division algorithm is adapted to calculate s-polynomial as defined in \ref{def2.7}. Division removes an appropriate multiple of $f$  from $g$ in order to cancel off $LT(f)$ with the term involving $h$ in $g$. We perform division as follows,
			$g-\frac{\lambda}{LC(f)}x_l f x_r=r$, where $\lambda \in K$ is to be chosen from $\{LC(x_l), LC(x_r)\}$. Similarly, for a set $F=\lbrace f_1, \dots, f_s \rbrace$; any $f \in A$ can be written as $f= u_{l1} f_1v_{r1} + \dots + u_{ls} f_sv_{rs} + r$ where $u_{li},v_{ri}, r, f_i \in A$ and either $r=0$ or $r$ is a linear combination, with coefficients in $K$, of monomials which are not divisible by $LT(f_i)$. $r$ is the remainder of $f$ after dividing by $F$. We denote by $r= Red_F(f)$ a reduction of $f$ with respect to the set $F$. Moreover if $u_{li}f_iv_{ri} \neq 0$ then we have that $LM(f) \succeq LM(u_{li}f_iv_{ri})$ of all $i$.
	
	\begin{definition}\label{def2.7} Let $f,g \in A$ and the leading monomials of $f$ and $g$ overlap such that $x_1LM(f)y_1=x_2LM(g)y_2$,where $x_1, x_2, y_1, y_2 \in M$ are chosen so that at least one of $x_1$ and $x_2$ and at least one of $y_1$ and $y_2$ is equal to unit monomial. Then the S-polynomial associated with this overlap is given by
		\[ S(f,g)=\lambda_1 x_1 \cdot f \cdot y_1 - \lambda_2 x_2 \cdot g \cdot y_2\]
		where $\lambda_1=\frac{LC(x_1)}{LC(f)}$ when $x_1 \neq 1$ or $\lambda_1=\frac{LC(y_1)}{LC(f)}$ when $y_1 \neq 1$ and  $\lambda_2=\frac{LC(x_2)}{LC(g)}$ when $x_2 \neq 1$ or $\lambda_2=\frac{LC(y_2)}{LC(g)}$ when $y_2 \neq 1$.
	\end{definition}
	
	\subsection{Mora's Algorithm}
	In commutative Gr\"{o}bner basis theory, Buchberger's Algorithm is used to compute the Gr\"{o}bner basis. Moreover, Dickson's Lemma and Hilbert's Basis Theorem assures termination of this Algorithm for all possible inputs. Our next result is Mora's Algorithm, that mimic Buchberger's Algorithm for noncommutative polynomial rings. However there is no analogous Dickson's Lemma for noncommutative monomial ideals, hence Mora's Algorithm does not terminates for all possible inputs.
	
	\begin{center}
		\begin{algorithm}[H]
			\caption{Noncommutative Mora's Algorithm \cite{mora}}
			\SetKwInOut{Input}{Input}
			\SetKwInOut{Output}{Output}
			\SetKwInOut{Initialize}{Initialize}
			\Input{A basis $F=\lbrace f_1, \dots, f_n \rbrace$ for ideal I over a non commutative polynomial ring $A=K[x_1, \dots, x_n]$ and an admissible order $\prec$.}
			\Output{A Gr\"{o}bner basis $G=\lbrace g_1, \dots, g_t \rbrace$ for I (In the case of termination)}
			{Let $G=F$ and let $B= \emptyset$.
				for each pair ($g_i, g_j) \in G, i \leq j)$ add an S-polynomial $s(g_i, g_j)$ to B for each overlap $x_1LM(g_i)y_1=x_2 LM(g_j)y_2$ between the leading monomials $LM(g_i)$ and $LM(g_j)$.}\\
			\While{$B\neq \emptyset$}{Remove the first entry $s_1$ from B. $s_1^{\prime}=Red_G(s_1)$
				
				\If{$s_1^{\prime} \neq 0$}{Add $s_1^{\prime} $ to G and then for all $g_i \in G$ add all $S(g_i,g_j)$ to B.}
			}
			{Return G}		
			
		\end{algorithm}
	\end{center}
	 The next result shows that it is indeed possible to have infinite Gr\"{o}bner basis for some finitely generated ideal of $A=K[x_1, \dots, x_n]$.
	\begin{proposition} Not all noncommutative monomial ideals are finitely generated.
		\proof Assume to the contrary that all noncommutative monomial ideals are finitely generated, and consider an ascending chain of such ideals $J_1 \subset J_2 \subset \dotsi .$ Then $J=\cup J_i$ is finitely generated and there is some $d \ge 1$ such that $J_d=J_{d+1} = \dotsi.$ For a counterexample, let $A =K[x, y]$	be a noncommutative polynomial ring, and define $J_i$ for ($i > 1$) to be the ideal in $A$	generated by the set of monomials $\lbrace xyx, xy^2x,\dots, xy^{i}x \rbrace$. Thus we have  an ascending chain of such ideals $J_1 \subset J_2 \subset \dotsi .$ However because no member of this set is a multiple of any other member of the set, it is clear that there cannot be a $d\geq 1$ such that $J_d=J_{d+1} = \dotsi$, because $xy^{d+1}x \in J_{d+1}$ and $xy^{d+1}x \notin J_d$ for all $d > 1$. $\square$
	\end{proposition}	
	
\subsection{Path Algebra}
	Now we are ready to study the core object of this paper, a noncommutative free associative algebras called path algebras. For a path algebra we will give concrete example of overlap relations and noncommutative Gr\"{o}bner basis. We begin by defining a path algebra.
	\begin{definition}A $quiver$ is a quadruple $Q=(Q_0, Q_1, s, t)$ consisting of two sets: $Q_0$ whose elements are called points or vertices, say $\lbrace 1, 2, 3,...., n, ..... \rbrace$,  $Q_1$ whose elements are called arrows, say $\lbrace \alpha_1, \alpha_2, \alpha_3, \ldots,\alpha_n, \ldots \rbrace$	and two maps: $s,t:Q_1 \longmapsto Q_0$ which associates to each arrow $\alpha \in Q_1$ its source $s(\alpha) \in Q_0$ and its target $t(\alpha) \in Q_0$ respectively.\\
	\end{definition}
	
	\begin{definition}
		An arrow $\alpha \in Q_1$ of source $s(\alpha)=1$ and target $t(\alpha)=2$ is usually denoted by $\alpha : 1 \longmapsto 2$. A $path$ $x$, of length $l >1$, with a source $a$ and target $b$, is a sequence of arrows $\alpha_1, \alpha_2, \alpha_3, ....... ,\alpha_n$ with $a= s(\alpha_1)$ and $b= t(\alpha_n$) where $ \alpha_k \in Q_1$ for all $1\leq k \leq n$ and $t(\alpha_k) = s(\alpha_{k+1})$ for $1 \leq k < n$. Such a path $x$ is denoted by $x= \alpha_1 \alpha_2 \alpha_3 .... \alpha_n$ and visualized as: \[ a= 1 \xrightarrow{\alpha_1} 2 \xrightarrow{\alpha_2} 3 \xrightarrow{\alpha_3} .... \xrightarrow{\alpha_n-1} n \xrightarrow{\alpha_{n}} n+1=b \] 
	\end{definition}
	
	\begin{definition}
		The $length$ of a path $x$, denoted by $l=l(x)$ is the number of arrows in it. An $arrow$ $\alpha : 1 \longmapsto 2$ is a path of length 1. A $trivial$ path denote by $v_i$ is a path of length zero associated with each vertex $i$. A path of length $l \geq 1$ is called a $cycle$ whenever its source and target coincide. A $loop$ is a cycle of $l= 1$. A quiver is said to be $acyclic$ if it has no cycles. A quiver is said to be $finite$ if $Q_0$ and $Q_1$ are both finite sets.
	\end{definition}
	
	\begin{definition} \label{2.7}Let $Q$ be a quiver and $K$ an arbitrary field. The path algebra $KQ$ of $Q$ is the $K$-algebra whose underlying $K$-vector space has as its basis the set of all paths of length $l \geq 0$ in $Q$, and such that the product of two basis vectors namely $x= \alpha_1 \alpha_2 \alpha_3 .... \alpha_n$ and $y= \beta_1 \beta_2 \beta_3 .... \beta_k$ is defined by 
		\[ 
		xy =
		\begin{cases}
		\alpha_1 \alpha_2 \alpha_3 ... \alpha_n \beta_1 \beta_2 \beta_3 ... \beta_k ,& \text{if} \quad s(y) =t(x) \\
		0,  \qquad \text{otherwise}
		
		\end{cases}
		\]
		i.e the product $xy$ is a concatenation or zero otherwise, so that $ Q \bigcup \lbrace 0 \rbrace$ is closed under multiplication. Multiplication as defined above is also distributive $K$-linearly in  $ Q \bigcup \lbrace 0 \rbrace$. Addition in $KQ$ is the usual $K$-vector space addition where $Q$ is a $K$-basis for $KQ$. 
	\end{definition} 
	The following two results, \ref{Remark2.13} and \ref{lemma2.14}, shows that $KQ$ as defined in \ref{2.7} is indeed an associative algebra.
	\begin{remark}{(Properties)}\label{Remark2.13}
		\begin{enumerate}[noitemsep]
			\item  Let $Q$ be finite. The set $\lbrace v_1, v_2, v_3, ....., v_n \rbrace$  of the trivial paths corresponding to the vertices $\lbrace 1, 2, \dots, n \rbrace$ is a complete set of primitive orthogonal idempotents. Thus $	  1= v_1+ v_2+ v_3+ .....+ v_n= \sum_{i=1}^{n} v_i $ is the called the identity element of $KQ$. 
			\item For each arrow $\alpha: 1 \longmapsto 2$ we have the following defining relations; 
			\begin{itemize}
				\item	  $ v_{i}^{2} = v_i v_i= v_i$ for $i=1,2.$
				\item $ v_1 \alpha = \alpha $
				\item $ \alpha v_2 = \alpha $
				\item $  v_1 v_2 = 0. $
			\end{itemize}
			\item Let $Q$ denote the set of all paths of length $l \geq 0$, then the above product extend to all elements of $KQ$ and there is a direct sum
			\[
			KQ=  KQ_1 \oplus KQ_2 \oplus \dots \oplus KQ_i \oplus \dots
			\]
			Where $KQ_i$ is subspace of $KQ$ generated by the set $Q_i$, where $Q_i$ is the set of all paths of length $i$,  over $K$. Since the product of path of length $n$ with path of length $m$ is zero or a path of length $n+m$ then the above decomposition defines a grading on $KQ$. Hence $KQ$ is a graded $K$-algebra.
		\end{enumerate}
	\end{remark}
	
	\begin{lemma}\cite{assem}\label{lemma2.14}
		Let $Q$ be a quiver and $KQ$ be its path algebra. Then
		\begin{enumerate}[noitemsep]
			\item $KQ$ is an associative algebra.
			\item $KQ$ has an identity element if and only if $Q$ is finite.
			\item $KQ$ is finite dimensional if and only if $Q$ is finite and acyclic.
		\end{enumerate}
	\end{lemma}
	\begin{definition} An element $f \in KQ;\quad (f=\sum \lambda_i x_i, \quad \lambda_i \in K)$, is a linear combination of paths $x_i \in Q$ over $K$. Elements of $KQ$ will be called polynomials. The paths $x_i \in Q$ appearing in each polynomials will be called monomials. We shall denote by $Mon(f)$ the set of all monomials $x_i$ appearing in the polynomial $f$.
	\end{definition}
	
	\begin{example}
			If $Q$ consist of one vertex and $n$ loops, $\alpha_1, \alpha_2 \dots  \alpha_n,$ then $KQ \cong K[X_1, X_2, \dots, X_n]$. 
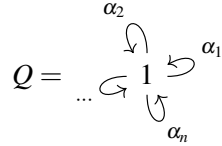
			\[Q=
			\begin{tikzcd}
			1 \arrow[out=0,in=30,loop,swap,"\alpha_1"]
			\arrow[out=90,in=120,loop,swap,"\alpha_2"]
			\arrow[out=180,in=210,loop,swap,"\cdots"]
			\arrow[out=270,in=300,loop,swap,"\alpha_n"]
			\end{tikzcd}
			\]
			The isomorphism is induced by the $K$-linear maps
			\[
			v_1 \longmapsto 1, \quad \alpha_1 \longmapsto X_1, \quad \alpha_2 \longmapsto X_2, \quad \dots \quad \alpha_n \longmapsto X_n.
			\]
		\end{example}

	\subsection{Basics to  Noncommutative Gr\"{o}bner Basis in a Path Algebra}
	For the rest of the paper, $Q$  will contain paths of finite length $l$. By convection we write a path $ \alpha_1 \alpha_2 \alpha_3 \dots \alpha_n $ from left to right such that $t(\alpha_i) = s(\alpha_{i+1})$. For path $x= \alpha_1 \alpha_2 \alpha_3 \dots \alpha_n$ we denote its length $l(x)=n$.
	\begin{definition}
		\begin{enumerate}[noitemsep]
			\item A subset $L$ of $KQ$ is called a left ideal if:
			\begin{enumerate}[noitemsep]
				\item  $0 \in L$
				\item $x +y \in L \quad \text{for all} \quad x,y \in L$ 
				\item $xy \in L \quad \text{for all} \quad x\in KQ \quad \text{and all} \quad y \in L$
			\end{enumerate}
			\item A subset $R$ of $KQ$ is called a right ideal if:
			\begin{enumerate}[noitemsep]
				\item  $0 \in R$
				\item $x +y \in R \quad \text{for all} \quad x,y \in R$ 
				\item $xy \in R \quad \text{for all} \quad x \in R \quad \text{and all} \quad y \in KQ$
			\end{enumerate}
			\item In general a subset $I$ of $KQ$ is called a two-sided ideal or simply an ideal, if it is both a left and a right ideal.
		\end{enumerate}
	\end{definition}
	In general an ideal $I$ in the path algebra $KQ$ has a Gr\"{o}bner basis depending on the ordering of the paths in $Q$.
	\begin{proposition}
		 An ideal $I$ in $KQ$ with some path ordering has a Gr\"{o}bner basis whenever the path ordering is admissible.
	\end{proposition}
	
	\begin{definition}[Path Ordering] By a path ordering we refer to the noncommutative ordering as defined in \ref{def2.2}. In addition we arbitrary order the vertices $ v_1 \prec v_2 \prec v_3 \prec \dots \prec v_k$ and arbitrary order the arrows all larger than given vertex say $v_k$ as $ v_k \prec \alpha_1 \prec \alpha_2 \prec \alpha_3 \prec \prec \alpha_r $.
	\end{definition}
	
	\begin{definition} \label{def534} A path order $\prec$ is said to be admissible order if:
		\begin{enumerate}[noitemsep]
			\item Whenever $x \neq y$ either $x \prec y$ or $x \succ y$.
			\item Every nonempty set of paths has a least element.
			\item $x \prec y \Rightarrow xz \prec yz$, whenever $xz \neq 0$ and $yz \neq 0$.
			\item Also $x \prec y \Rightarrow wx \prec wy$, whenever $wx \neq 0$ and $ wy \neq 0$.
			\item $x = yz$ implies $x \succeq y$ and $ x \succeq z$.
		\end{enumerate}
	\end{definition}
	
	\begin{remark} 
		Conditions $a$ through $c$ makes $\prec$ $a$ right admissible order. Condition $a, b$ and $d$ make $\prec$ left admissible ordering whilst condition $b$ say that an admissible ordering is a well ordering.
	\end{remark}
	
	While in section \ref{sec2.1} we defined a noncommutative monomial ordering, those monomial ordering were not necessarily admissible. In the following examples we use appropriate monomial (path) ordering to construct admissible ordering for paths in $KQ$.
	
	\subsubsection{Constructing Admissible Path Ordering}
		\begin{enumerate}[noitemsep]
			
			\item Left Lexicographic order: Let $x=\alpha_1 \dots \alpha_n$ and $y= \beta_1 \dots\beta_m$ be paths. We say that $x$ is less than $y$ with respect to left lexicographic order and denote $x \prec_{llex} y$ if there exist a path $z$ (otherwise we set $z=1$), such that $x=z \alpha_k \dots \alpha_n, \quad y=z\beta_s \dots \beta_m$ and $ \alpha_k \prec \beta_s$. Left lexicographic order is not a left admissible ordering since it is not a well ordering. For example let $Q$ be
			\[Q=
			\begin{tikzcd}
			1 \arrow[out=90,in=180,loop, swap, "\alpha"] \arrow[r, "\beta"] & 2
			\end{tikzcd}
			\]
			with $\alpha \prec \beta$. We have $( \alpha \beta \succ_{llex} \alpha^{2} \beta \succ_{llex} \alpha^{3} \beta \dots)$. Then the subset $\{  \alpha^{n} \beta : n \in \mathbb{N} - \lbrace 0 \rbrace \} \subset Q$ does not have a least element.
			
			\item Length  left lexicographic order: Let $x=\alpha_1 \dots \alpha_n$ and $y= \beta_1 \dots\beta_m$ be paths. We say that $x$ is less than $y$ with respect to length left lexicographic order and denote $x \prec_{Lex} y$ if $l(x) < l(y)$ or $l(x)= l(y)$ and $x \prec_{llex} y$. Length Left lexicographic order is a left admissible order.
			
			\item Right lexicographic order:   Let $x=\alpha_1 \dots \alpha_n$ and $y= \beta_1 \dots\beta_m$ be two paths in $Q$. We say that $x$ is less than $y$ with respect to right lexicographic order and denote $x \prec_{rlex} y$ if there exist a path $z$ (otherwise we set $z=1$), such that $x= \alpha_1 \dots \alpha_k z, \quad y=\beta_1 \dots \beta_s z$ and $ \alpha_k \prec \beta_s$.
			This ordering is not a well ordering and hence not admissible.
			
			\item Length right lexicographic order: Let $x=\alpha_1 \dots \alpha_n$ and $y= \beta_1 \dots\beta_m$ be paths. We say that $x$ is less than $y$ with respect to length right lexicographic order and denote $x \prec_{rLex} y$ if $l(x) < l(y)$ or $l(x)= l(y)$ and $x \prec_{rlex} y$. Length right lexicographic order is a right admissible order.
	\end{enumerate}
	By length lexicographic ordering, we mean an ordering which is both Length light lexicographic order and Length right lexicographic order. With admissible ordering we are glad to attempt calculating Gr\"{o}bner basis in path algebras. Calculating this basis is a series of division and reduction as summarized below.
	
	\begin{definition} Let $\prec$ be an admissible ordering and $A=KQ$. Then definition \ref{def2.5} hold true for all $f \in A$. Moreover $x$ left divide $y$ if $y=wx$ for some path $w\in Q$. $x$ right divide $y$ if $y=xz$ for some path $z \in Q$. Hence  $x$ divides $y$ if $y=wxz$ for some paths $w,z \in Q$.
	\end{definition}
	
	\begin{definition}Let $x$ and $y$ be paths. 
	 An element $f \in KQ \setminus \lbrace 0 \rbrace$ is said to be an uniform if there exist vertices $u$ and $v$ such that $f=uf = fv =ufv$.
		
	\end{definition}
	
	\begin{proposition} \cite{micah} All elements of $KQ$ are uniform. 
		\proof $f= \sum_{i=1}^{n} \lambda_i x_i$ is uniform since for each monomial $x_i$, which is a sequence of arrows, has a source vertex say $u_i$ and a target vertex say $v_i$ and hence $x_i = u_i x_i v_i$. Therefore $f$ is sum if uniform elements $f = \sum\limits_{i,j=1}^{n} u_i f v_j$. $\square$
	\end{proposition}
	
	\begin{definition} Let $H$ be a subset of $KQ$ and $g \in KQ$. We say that $g$ can be reduced by $H$ if for some $x \in$ Mon$(g)$ there exist $h \in H$ such that $LM(h)$ divides $x$, i.e $x= p LM(h) q$ for some monomials $p,q \in KQ$. The reduction of $g$ by $H$ is given by $g- \lambda p h q$ where $h \in H$, $p,q \in Q$ and $\lambda \in K \setminus \lbrace 0 \rbrace$ such that $\lambda p LM(h) q$ is a term in $g$, $\lambda$ is uniquely determined by $\lambda =\frac{LC(g)}{LC(h)}.$ Moreover;
		\begin{enumerate}
			\item A total reduction of $g$ by $H$ is an element resulting from a sequence of reductions that cannot be further reduced by $H$.
			\item We say that an element $g \in KQ$ reduces to $0$ by $H$ if there is a total reduction of $g$ by $H$ which is $0$. In general two total reductions need not be the same.
			\item A set $H \subset KQ$ is said to be a reduced set if for all $g \in H$ , $g$ cannot be reduced by $H- \lbrace g \rbrace$.
		\end{enumerate}
	\end{definition}

	\section{One-side Gr\"{o}bner Bases in Path Algebra}
	We now look at the left and right division algorithms in a path algebra. These algorithms will be an indirect entries in the respective left and right Buchberger's Algorithm which in-turn produces respective left and right Gr\"{o}bner basis. Onesided Gr\"{o}bner basis and all one-sided algorithms are given in this section.
	
	\subsection{Left Gr\"{o}bner Bases in Path Algebra}
	
	\begin{definition} Let $L$ be a left ideal of $KQ$ with a left admissible order $\prec$. We say that a set $G_L \subset  L$ is a left Gr\"{o}bner basis for $L$, if for all $f \in$ L $\setminus \lbrace 0 \rbrace$ there exist $g \in G_L$ such that $LM(g)$ left divides $LM(f)$. Equivalently we say that a set $G_L \subset  L$ is a left Gr\"{o}bner basis for $L$ with respect to a left admissible order $\prec$ if $\langle LM(G_L) \rangle = \langle LM(L)\rangle$.
	\end{definition}
	
	\begin{theorem}\label{thm612} \cite{attan} Let $\prec$ be a left admissible ordering and $S = \lbrace  f_1, \dots, f_n \rbrace$ be a set of nonzero polynomials in $KQ$. For $g \in KQ \setminus \lbrace 0 \rbrace$ there exist a unique determined expression $g=\sum\limits_{i=1}^{n} g_i f_i + h$ where $h, g_1, \ldots, g_n \in KQ$ satisfying:
		\begin{itemize}
			\item[A.] For any path $p$ occurring in each $g_i, \quad t(p)= s(LM(f_i)).$
			\item[B.] For $i>j$, no term $g_i LT(f_i)$ is left divisible by $LT(f_j).$
			\item[C.] No path in $h$ is left divisible by $LM(f_i)$ for all $1 \leq i \leq n$.
		\end{itemize}
	\end{theorem}
	\begin{remark} The expression $g=\sum\limits_{i=1}^{n} g_i f_i + h$ in  theorem \ref{thm612} is called $left$ $standard$ representation of $g \in KQ$ with respect to the set $S$. Algorithm \ref{alg05} gives as an output $h$, a remainder of $g$ after left division by $S$. We denote by $LRed_S(g) = h$ the particular remainder of $g$ by the set S produced by division algorithm with respect to a fixed admissible ordering.
	\end{remark}
	
	\begin{center}
		\begin{algorithm}[H]
			\caption{Left Division Algorithm}
			\label{alg05}
			\SetKwInOut{Input}{Input}
			\SetKwInOut{Output}{Output}
			\SetKwInOut{Initialize}{Initialize}
			\Input{$g, S=\lbrace f_1, \ldots, f_n \rbrace$ $g,f_i \in KQ\setminus  \lbrace 0 \rbrace$ and left admissible order $\prec$ on $KQ$}
			\Output{$g_i, \ldots, g_n, h \in KQ$ such that $g=\sum\limits_{i=1}^{n} g_i f_i + h$ }
			\begin{enumerate}[noitemsep]
				\item[a.] {For any multiple $v_{1i}$ of $LT(f_1)$ occurring in $g$ with $( 1 \leq i \leq r_1)$, find for each $i$ a term $h_{1i}$ such that $v_{1i} = h_{1i} LT(f_1)$. Afterwards do the same for any multiple $v_{2i}$ of $LT(f_2)$ occurring in $g$ such that $v_{2i} = h_{2i} LT(f_2)$ with $1\leq i\leq r_2$. Continue in this way for any multiple $v_{ki}$ of $LT(f_k)$ such that $v_{ki} = h_{ki}LT(f_k)$ with $1 \leq i \leq r_k$ and $k \in {3, \ldots, n}$} \\
				\item[b.] {Write $g = \sum\limits_{j=1}^{n}(\sum\limits_{i=1}^{r_j}h_{ji})LT(f_j) + h_1$ and set $g^1 = g -(\sum\limits_{j=1}^{n}(\sum\limits_{i=1}^{r_j}h_{ji})f_j + h_1)$} \\
				{
					\item[c.] If $g^1 = 0$ then we are done and $g = \sum\limits_{j=1}^{n} g_j LT(f_j) +h_1$ where $	g_j = \sum\limits_{i=1}^{r_j}h_{ji}$ and $h_1 = h$.
				}\\
				{
					\item[d.] If $g^1 \neq 0$, go back to $a$ and continue the process with $g=g^1$	
				}
			\end{enumerate}
		\end{algorithm}
	\end{center}
	
	\begin{example} 
		Let $Q$ be the quiver with one vertex and three loops over the field of rationals.
		\[Q=
		\begin{tikzcd}
		1 \arrow[out=0,in=30,loop,swap,"x"]
		\arrow[out=90,in=120,loop,swap,"y"]
		\arrow[out=180,in=210,loop,swap,"z"]
		\end{tikzcd}
		\]
		With a left length lexicographic ordering $z\prec y \prec x$. We find the standard representation of $g=zxxyz + xyxxy - xyz$ with respect to the set $\lbrace f_1=xyz-zy, f_2=xxy-yx \rbrace$. We first note that $LM(f_1)=xyz$ and $LM(f_2)=xxy.$ Initializing we get $g=zxLM(f_1) + xyLM(f_2)-LM(f_1).$ We replace $g$ by $g^1=g-(zxf_1 +xyf_2- f_1)=zxzy +xyyx +zy$. Neither $LM(f_1)$ and $Lm(f_2)$ left divides $zxzy +xyyx +zy$, so we set $h=zxzy +xyyx +zy$ and $zxzy +xyyx +zy$ is replaced by $0$ and the algorithm stops. Thus the standard representation of $g$ is  $g=zxf_1 +xyf_2- f_1 + h$.
	\end{example}
	
	\paragraph{Proof of Theorem \ref{thm612} }
	\begin{enumerate}[noitemsep]
		\item \underline{Existence:} First the algorithm removes any multiple of $f_1$ from $g$. Then removes any multiple of $f_2$ and continue in this way until any multiple of any of $f_k$ has been removed. In this case if $g=\sum_{j=1}^{n} \sum_{i=1}^{r_j}h_{ji}LT(f_i) + h_1$ is the resulting standard representation of $g$, we have either $g^1= g-(\sum_{j=1}^{n} \sum_{i=1}^{r_j}h_{ji}(f_i) + h_1) = 0$ Or $LM(g^1) \prec LM(g)$. Since the path ordering $\prec$ is well ordering, by recursion the algorithm produces a standard representation for $g^1$, $g^1= \sum_{j=1}^{n} \sum_{i=1}^{r_j}h_{ji}^{1}(f_i) + h^1$ satisfying conditions $A$, $B$ and $C$. Thus $g=\sum_{j=1}^{n} \sum_{i=1}^{r_j}(h_{ji}+ h_{ji}^{1}) (f_i) + (h_1 + h^1)$ is a representation for $g$ satisfying the conditions $A$, $B$ and $C$. \item \underline{Uniqueness:} For $g \in L \setminus \lbrace 0 \rbrace$ let $g= g_1 f_1 + \dots+ g_n f_n + h$. Then the three conditions $A$, $B$ and $C$ implies that the terms $LT(g_if_i)=LT(g_i)LT(f_i)$ and $LT(h)$ do not divide each other to the left. Otherwise these terms cancels with each other into zero polynomial. Therefore the representation $g=\sum_{i=1}^{n}g_if_i + h$ is unique. 
		\item \underline{Termination:} The algorithm produces elements $g, g^1, g^2, \dots, g^k$ so that at each $k^{th}$ iteration $LM(g^{k+1}) \prec LM(g^{k})$. Since $\prec$ is a well ordering, the algorithm terminates at some $g^k = 0$ satisfying the conditions of the theorem. $\square$
	\end{enumerate}

	Given a finite generating set $S = \{f_1,\ldots, f_n\}$. For a left admissible order $\prec$, the following algorithm gives as an output $R_L=R_L(S)$, a left reduction of $S$.

	\begin{center}
		\begin{algorithm}[H]
			\caption{Set Left Reduction Algorithm}
			\SetKwInOut{Input}{Input}
			\SetKwInOut{Output}{Output}
			\SetKwInOut{Initialize}{Initialize}
			\Input{$S= \lbrace f_1, \ldots, f_n \rbrace$ $f_i \neq 0$ and a left admissible ordering $\prec$}
			\Output{$R_L$ a left reduction of the set $S$}
			\begin{enumerate}[noitemsep]
				\item[a.] $R_L=\emptyset$\\
				\item[b.] {Find the maximal element $f_k$ of $S$ with respect to $\prec$, for $1 \leq k \leq n.$}\\
				\item[c.] {Write $S=S-\lbrace f_k \rbrace$} \\
				\item[d.] {Do $f_{k}^{\prime}= LRed_{S \cup R_L}(f_k)$}\\
				\item[e.] {If $f_{k}^{\prime}\neq 0$ then $R_L=R_L \cup \lbrace \frac{f_{k}^{\prime}}{LM(f_{k}^{\prime})} \rbrace$ }\\
				\item[f.] { If $f_k \neq f_{k}^{\prime}$ then $S=S \cup R_L$}; Go back to $a$ and continue with the process. 
			\end{enumerate}
			
		\end{algorithm}
	\end{center}
	
	\begin{proposition} \cite{attan} Let $G=\lbrace f_1, \dots, f_n \rbrace \subset KQ$ be a left Gr\"{o}bner basis for the ideal $L=\langle f_1, \dots, f_n \rangle \subset KQ$. If $g=\sum \limits_{i=1}^{n}g_i f_i + h$ is a left standard expression of $g \in KQ \setminus \lbrace 0 \rbrace$ then $g \in L$ if and only if $h=0$.
		\proof If $h=0$ clearly $g \in L$. Conversely if $g \in L$ then $h \in L \Longrightarrow LM(h) \in \langle LM(f_1), \dotsi, LM(f_n) \rangle$ which is impossible by the theorem \ref{thm612}.
	\end{proposition}
	
	\begin{definition}[Left S-Polynomial] Let $f,g \in KQ \setminus \lbrace 0 \rbrace$ and $\prec$ be a left admissible ordering. Let $p, q$ be paths such that $p LM(f)= q LM(g)$, the left S-polynomial $S_L(f,g)$ is defined as
		\[ S_L(f,g) = \frac{p}{LC(f)}\cdot f - \frac{q}{LC(g)}\cdot g
		\]
	\end{definition}
	
	\begin{theorem}[Left Buchberger's Criterion]\cite{attan} Let $f_1, \dots, f_n \in KQ \setminus\lbrace 0 \rbrace$ and $\prec$ be a left admissible ordering. Let $S_L(f_i,f_j)= \sum_{k=1}^{n}g_k f_k +h_{ij}$ be a left a standard expression of $S_L(f_i,f_j)$ for each pair $(i,j)$.  $\{f_1, \dots, f_n\}$ form a left Gr\"{o}bner basis for $L=\langle  f_1, \dots, f_n \rangle$ if and only if all the remainders $h_{ij}$ are zero.
	\end{theorem}
	
	\begin{center}
		\begin{algorithm}[H]
			\caption{Left Buchberger's Algorithm}
			\label{alg07}
			\SetKwInOut{Input}{Input}
			\SetKwInOut{Output}{Output}
			\Input{$L = \langle f_1, \dots, f_n \rangle \subset KQ$ and a left admissible order $\prec$.}
			\Output{A reduced left Gr\"{o}bner basis $G_m$ for $L$.}
			\begin{enumerate}[noitemsep]
				\item[a.] {$m=0$; $G_0=\emptyset$; $G_1= R_L( \lbrace f_1, \dots, f_n \rbrace)$}\\
				\item[b.] {While $G_m \neq G_{m+1},$ $m=m+1$}\\
				\item[c.] {For all $g,h \in G_m$ find all $S_L(g,h) \neq 0$}\\
				\item[d.] {Write $G_{m}^{\prime} =G_{m}^{\prime} \cup \lbrace S_L(g,h)\rbrace $}\\
				\item[e.] {$G_{m+1}=R_L(G_{m}^{\prime})$}
			\end{enumerate}
		\end{algorithm}
	\end{center}
	
	\subsection{Right Gr\"{o}bner Basis in a Path Algebra}
	The right division and reduction algorithms in many ways gives a "right" version of the Gr\"{o}bner basis theory as discussed in the previous section. This means that concepts from the previous section  will have to be duplicated with slight variant in the division as right sided operation. It is hence omitted here and we instead illustrate right Gr\"{o}bner basis in a path algebra using the following example.
	\begin{example} 
		Let $Q$ be the quiver; $\begin{tikzcd}
		1 \arrow[r,"z"] \arrow[dr,"y"] & 2 \arrow[d,"t"] \\ & 3 \arrow[out=60,in=330,loop,"x"]
		\end{tikzcd}$\\
		Let $F = \lbrace f_1 = ztx^3, f_2 = zt + y \rbrace$ be a subset of $KQ$ with respect to the right length lexicographic ordering $v_1 \prec v_2 \prec v_3 \prec t \prec z \prec y \prec x$. Running Algorithm \ref{alg10}, we note that $LM(f_1)=ztx^3$ and $LM(f_2)=zt$ and they only factor each other to the right in one way namely $LM(f_1)v_3 = LM(f_2)x^3$. Thus we have one right S-polynomial $S_R(f_1,f_2)=f_1 v_3 -f_2 x^3= -yx^3$. Neither $LM(f_1)$ nor $LM(f_2)$ right divide $-yx^3$ so we add $f_3=-yx^3$ to F. Now every right S-polynomial reduces to zero by F. Thus $F=\lbrace f_1, f_2, f_3 \rbrace$ is a right Gr\"{o}bner basis for the ideal $R=\langle f_1, f_2 \rangle$.
	\end{example}
\section{Twosided Gr\"{o}bner Bases}

We say that a set $G \subset I$ is a Gr\"{o}bner basis for $I$ with respect to an order $\prec$ if $\langle LM(G)\rangle = \langle LM(I)\rangle$.

\begin{proposition}
	If $G$ is a Gr\"{o}bner basis for the ideal $I$, then $G$ is a generating set for the elements of $I$ and also $G$ reduces elements of $I$ to $0$.
	\proof Let $KQ$ be a path algebra with an admissible ordering $\prec$. Let $I$ be an ideal and let $G$ be a Gr\"{o}bner basis for $I$. Let $f_{i} \in I$, $i=1, \ldots,n,\ldots $, for every $f_n \in I$ such that $f_n \neq 0 \quad \exists \quad g \in G$ such that $LM(g)$ divides $LM(f_n)$. Let $f_{n+1} = f_n - \frac{LC(f_n)}{LC(g)} x g y$ be a reduction of $f_n$ by $g$. Then $LM(f_{n+1}) \prec LM(f_n)$. But $g, f_n \in I \quad \Longrightarrow f_{n+1} \in I$. Repeating this reduction on $f_i$ to produce $f_{i+1}$ yields a decreasing sequence $
	LM(f_1) \succ LM(f_2) \succ \dots \dots, $ which terminates only if $f_n = 0$. Since $\prec$ is an admissible order, every set of paths has a least element hence the sequence must terminate with $f_n = 0.$ $\square$
\end{proposition}

\subsection{Division Algorithms}
\begin{theorem} \label{thm711} Let $\prec$ be an admissible ordering and $S = \lbrace  f_1, \ldots, f_n \rbrace$ be a set of non zero polynomials in $KQ$. For $g \in KQ \setminus \lbrace 0 \rbrace$ there exist a unique determined expression $g=\sum\limits_{i=1}^{n} w_i f_i z_i + h$ where $h, w_1, \ldots, w_n, z_1 \ldots, z_n \in KQ$ satisfying:
	\begin{itemize}
		\item[A3.] For any path $p$ occurring in each $w_i$, $ t(p)= s(LM(f_i))$ and for any path $q$ occurring in $z_i,$  $t(LT(f_i))=s(q).$ 
		\item[B3.] For $i>j$ no term $w_i LT(f_i) z_i$ is divisible by $LT(f_j).$
		\item[C3.] No path in $h$ is divisible by $LM(f_i)$ for all $1 \leq i \leq n$.
	\end{itemize}
\end{theorem}

\begin{center}
	\begin{algorithm}[H]
		\caption{Twosided Division Algorithm}
		\label{alg11}
		\SetKwInOut{Input}{Input}
		\SetKwInOut{Output}{Output}
		\SetKwInOut{Initialize}{Initialize}
		\Input{$g, S=\lbrace f_1,..., f_n \rbrace$ and an admissible order $\prec$ on elements of $KQ$.}
		\Output{$w_1, ..., w_n, z_1, ..., z_n, h \in KQ$ such that $g = \sum\limits_{i=1}^{n} w_i f_i z_i + h.$}
		\begin{enumerate}[noitemsep]
			\item[a.]	{For any multiple $O_{1i}$ of $LT(f_i)$ occurring in $g$ with $1 \leq i \leq r_1,$ find for each $i$ the terms $u_{1i}$ and $v_{1i}$ such that $O_{1i}=u_{1i} LT(f_1)v_{1i}$. Following this do the same for any multiple $O_{2i}$ of $LT(f_2)$ occurring in $g$ such that $O_{2i}=u_{2i} LT(f_2)v_{2i}$ with $1 \leq i \leq r_2.$ Continue in this way for any multiple $O_{ki}$ of $f_k$ such that $O_{ki}=u_{ki} LT(f_k)v_{ki}$ with $1 \leq i \leq r_k$ and $k \in \lbrace 3, ..., n \rbrace.$}\\
			\item[b.]	{Write $g=\sum_{j=1}^{n} \sum_{i=1}^{r_j}u_{ji}LT(f_j)v_{ji} + h_1$ and set $g^1= g -(\sum_{j=1}^{n} \sum_{i=1}^{r_j}u_{ji}f_j v_{ji} + h_1)$}\\
			\item[c.]	{If $g^1= 0$ then we are done and $g=\sum\limits_{j=1}^{n} w_j f_j z_j + h_1$ where $w_j=\sum_{i=1}^{r_j} u_{ji}$, $z_j=\sum_{i=1}^{r_j }v_{ji}$ and $h=h_1$}\\
			\item[d.] 	{If $g^1\neq 0,$ go back to $a$ and proceed with $g= g^1.$}
		\end{enumerate}	
	\end{algorithm}
\end{center}
Let $Red_S(g)=h$ denote the particular total reduction of an element $g$ by a set $S$ produced by the algorithm \ref{alg11} with respect to a fixed admissible ordering.
\paragraph{Proof}
\begin{enumerate}[noitemsep]
	\item \underline{$existence:$} This algorithm finds a standard representation of $g$ as follows. First it removes any multiple of $f_1$ in $g$. Afterwards removes any multiples of $f_2$. Continue in this way until any multiple of any $f_k, k \in \{ 3, 4, \dots, n \}$ has been removed. Hence if $g=\sum_{j=1}^{n} \sum_{i=1}^{r_j}u_{ji}LT(f_i)v_{ji} + h_1$ is the resulting representation of $g$ then either $g^1=g-(\sum_{j=1}^{n} \sum_{i=1}^{r_j}u_{ji}(f_i)v_{ji} + h_1)$ equal to zero and we are done, or $LM(g)\succ LM(g^1)$. Since $\prec$ is a well ordering then the algorithm finds a representation $g^1= \sum_{j=1}^{n} \sum_{i=1}^{r_j}u_{ji}^{1} f_i v_{ji}^{1} + h^{1}$ satisfying conditions $A3$, $B3$ and $C3$ so that $g= \sum_{j=1}^{n} \sum_{i=1}^{r_j}(u_{ji}^{1} + u_{ji}) f_i (v_{ji}^{1}+v_{ji}) + (h^{1}+ h_1)$ is the standard representation of $g$ satisfying conditions $A3, B3$ and $C3$.
	\item \underline{$Uniqueness:$} Given $g$ and conditions $A3$, $B3$ and $C3$, no term $LT(w_i f_i z_i)$  for all $1 \leq i \leq n$ divides $LT(h)$. Therefore the algorithm produces a unique standard representation $g=\sum_{i=1}^{n}w_i f_i z_i + h$ where $w_i$ or $z_i$ may be unit monomials.
	\item \underline{$termination:$} Note that the algorithm produces elements $g, g^1, g^2, ...., g^k$ such that at each $k^{th}$ iteration $LM(g^k) \succ LM(g^{k+1})$ and the algorithm must terminate at some $k$ where $g^k=\sum_{j=1}^{n} \sum_{i=1}^{r_j}u_{ji}(f_i)v_{ji} + h_k=0$ and every monomial occurring in the final $h_k$ is  not divisible by $LM(f_i), 1 \leq i \leq n$.		$\square$
	
\end{enumerate}

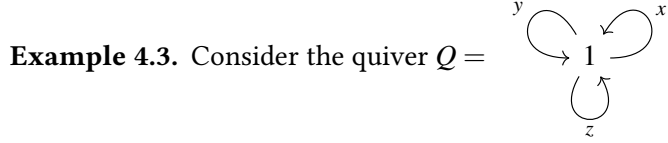
\begin{example} Consider the quiver $Q=
	\begin{tikzcd}
	1 \arrow[out=0,in=60,loop,swap,"x"]
	\arrow[out=120,in=180,loop,swap,"y"]
	\arrow[out=240,in=300,loop,swap,"z"]
	\end{tikzcd}
$\\
	
	Let $\prec$ the	length lexicographic order with $x \succ y \succ z$. Let's divide $f_1 = xy-x, f_2 = xx-xz$ into $f= zxxyx$. Note that the $LM(f_1)= xy$ and $LM(f_2)=xx$. Beginning the algorithm \ref{alg11}, we see that $zxxyx = (zx)LM(f_1)(x).$ Thus $p_1 =zx, q_1= x$ and we replace $zxxyx$ by $zxxyx - zx(f_1)x = zxxx$.	Now $LM(f_1)$ does not divide $zxxx$. Continuing, $LM(f_2)$ does. There are two ways to divide $zxxx$ by $xx$ and for the algorithm to be precise we must choose one. Say we choose the "left most" division. Then $zxxx = z(LM(f_2))x $ and we let $p_2= z, q_2 = x$ and replace $zxxx$ by $zxxx - z(f_2)x = zxzx$. Neither $LM(f_1)$ nor $LM(f_2)$ divide $zxzx$ so we let $r = zxzx$ and $zxzx$ is replaced by $0$ and the algorithm stops. We have $zxxyx = (zx)f_1 (x) + (z)f_2 (x) + zxzx$. The remainder is $zxzx$.
\end{example}

Given a finite generating set $S= \lbrace f_1, \ldots , f_n \rbrace$. For an ideal $I \subset KQ$ and an admissible order $\prec$, the following algorithm gives as an output $R(S)$ a finite monic reduced generating set for $I$.

\begin{center}
	\begin{algorithm}[H]
		\caption{Set Reduction Algorithm}
		\label{alg12}
		\SetKwInOut{Input}{Input}
		\SetKwInOut{Output}{Output}
		\SetKwInOut{Initialize}{Initialize}
		\Input{$S=\lbrace f_1, \ldots, f_n \rbrace$  $f_i \neq 0$ and an admissible ordering $\prec$.}
		\Output{$R=R(S)$ a reduction of elements of $S$}
		\begin{enumerate}[noitemsep]
			\item[a.] { $R=\emptyset $}\\
			\item[b.] {Find the maximal element $f_k$ of $S$ with respect to $\prec$.}\\
			\item[c.] {Write $S=S-\lbrace f_k \rbrace$}\\
			\item[d.] {Do $f_{k}^{\prime}=Red_{S \cup R}(f_k)$}\\
			\item[e.] {If $f_{k}^{\prime} \neq 0$ then $R=R \bigcup \lbrace \frac{f_{k}^{\prime}}{LM(f_{k}^{\prime})} \rbrace$}\\
			\item[f.] {If $f_k \neq f_{k}^{\prime}$ then $S= S \cup R$; Go back to $a$ and continue with the process.}
		\end{enumerate}
	\end{algorithm}
\end{center}

\begin{proposition} Given an ideal $I$ in $KQ$ and admissible order $\prec$, there is a unique Gr\"{o}bner basis $G$ such that $G$ is a reduced set and the coefficient of the leading monomials of the polynomials in $G$ are all $1$.
	\proof
	Let $KQ$ be a path algebra, $I$ an ideal and $\prec$ an admissible order. Let $G$ and $G^{\prime}$ be Gr\"{o}bner bases for $I$. Suppose  $G$ and $G^{\prime}$ are both reduced monic sets. Since $G \subset I$, for every $g_1 \in  G$ there exist $g^{\prime} \in G^{\prime}$ such that $LM(g^{\prime})$ divides $LM(g_1)$. Also since $G^{\prime}\subset I$ there exist $g_2 \in G$ such that $LM(g_2)$ divides $LM(g^{\prime})$. Thus $LM(g_2)$ divides $LM(g_1)$. But $G$ is a reduced set hence we must have that $g_ 2 = g_1$ so that $LM(g_1) = LM(g^{\prime}) = LM(g_2)$. So there is a bijection correspondence between elements of $G$ and the elements of $G^{\prime}$ with the same leading monomials. Thus $g^{\prime}$ cannot be reduced by $G - \lbrace g_1 \rbrace$. Hence $g^{\prime} - g_1$ cannot be reduced by G, since  $g^{\prime} - g_1 \in I$. Thus $g^{\prime} - g_1 = 0 \Longrightarrow g^{\prime} = g_1$ hence $G^{\prime}= G$. $\square$
\end{proposition}		

We call the unique reduced monic Gr\"{o}bner basis the reduced Gr\"{o}bner basis for $I$. The reduced Gr\"{o}bner basis $G$ is minimal in the sense that for any other reduced Gr\"{o}bner basis $G^{\prime}$ for the same ideal with the same admissible order, we have $LM(G^{\prime}) \subset LM(G)$.
\subsection{Twosided S-Polynomial}
While noncommutative S-polynomials for each pair of polynomials $f, g \in KQ$ may be different due to different factorizations in the set of monomials in $KQ$, for onesided case these S-polynomials are finitely many. However we may have ambiguity while dealing with twosided S-polynomials due to possible different choices of right and left factors of each overlap of $LM(f)$ and $LM(g)$. Therefore a condition namely $l(p) \leq l(LM(g))$ whenever $LM(f)\cdot p = q \cdot LM(g)$, is added to the definition of two sided S-polynomial to eliminate such ambiguity.

\begin{definition} Let $f, g \in KQ$ with an admissible order $\prec$ on elements of $KQ$. An $(f-g)$ overlap is said to occur if there are paths $p$ and $q$ of positive length such that $LM(f)p = q LM(g)$ where $l(p) \leq l(LM(g))$. Thus an $f$ and $g$ are said to have an overlap relation or a twosided S-polynomial denoted by $S(f,g)$ and defined as;
	\[
	S(f,g, p,q)=\frac{1}{LC(f)} f \cdot p - \frac{1}{LC(g)}q \cdot g.
	\]
\end{definition}

\begin{remark} Given elements $f, g \in KQ$ such that $LM(f)p = q LM(g)$ where $l(p) \leq l(LM(g))$, monomials $p$ and $q$ will not necessarily be unique. Consequently the same two elements $f$ and $g$ may still have multiple S-polynomials. In addition an element may have an S-polynomial with itself, i.e $S(f,f)$ will be a possible.
\end{remark}

\begin{example}

Let $Q$ be 
		$
		\begin{tikzcd}
		1 \arrow[out=0,in=90,loop,swap,"x"]
		\arrow[out=270,in=180,loop,"y"]
		\end{tikzcd}
		$ and $x \prec y$ with respect to the length lexicographic order. Let $f = 5yyxyx - 2xx$ and  $g = xyxy - 7y$. We see that $LM(f) = yyxyx$ and $LM(g)=xyxy$. The following are the S-polynomials among $f$ and $g$ are:
			\[\begin{array}{ccc}
			S(f,g, y, yy) =& \frac{1}{5}fy-yyg =& -\frac{2}{5}xxy + 7yyy	\\
			S(f,g,yxy,yyxy)=&\frac{1}{5}fy-yyg =&-\frac{2}{5}xxyxy + 7yyxyy\\
			S(g,g,xy,xy) =& gxy - xyg =& -7yxy + 7xyy 
						\end{array}\]
	\end{example}

\begin{lemma}[Bergman's Diamond, \cite{bergman}] \label{lemma4.8} Let $G$ be a set of uniform elements that form a generating set for the ideal $I\subset KQ$, such that for all $g,g_1 \in G$, $LM(g) \not | LM(g_1).$ If for each $f \in I$ and $g \in G$ every S-polynomial $S(f,g,p,q)$ is reduced to $0$ by $G$, then $G$ is a Gr\"{o}bner basis for $I$.
\end{lemma}
The beauty of algorithm \ref{alg11} in \ref{thm711} is that its outputs are uniform elements of $KQ$. In the next section we shall use \ref{lemma4.8}  to greatly reduces calculations and ascertain finite Gr\"{o}bner basis whenever algorithm \ref{alg13} terminates.
\subsection{The Main Theorem}

\begin{theorem}
	\label{thm731} Given a path algebra $KQ$, an admissible order $\prec$ and a finite generating set $\lbrace f_1, f_2, \dots , f_m \rbrace$ for an ideal $I$ the following algorithm gives a reduced Gr\"{o}bner basis for $I$ in the limit.
\end{theorem}

\begin{center}
	\begin{algorithm}[H]
		\caption{Twosided Buchberger's Algorithm}
		\label{alg13}
		\SetKwInOut{Input}{Input}
		\SetKwInOut{Output}{Output}
		\SetKwInOut{Initialize}{Initialize}
		\Input{$I=\langle f_1, \dots , f_n \rangle$, $f_i \neq 0$ and an admissible order $\prec$.}
		\Output{A reduced Gr\"{o}bner basis $G_m$ for $I$.}
		\begin{enumerate}[noitemsep]
			\item[a.] { $m=0; G_0 = \emptyset; G_1= R(\lbrace f_1, f_2, \dots , f_n \rbrace)$}\\
			\item[b.] {For $G_m \neq G_{m+1}; m=m+1$}\\
			\item {For all pairs $(g_i,g_j) \in G_m$ and all $1 \leq i \leq j \leq n$, find $S(g_i,g_j,p,q) \neq 0$}\\
			\item[c.]  {Do $G_{m}^{\prime}= G_m \bigcup \lbrace S(g_i,g_j,p,q) \rbrace$}\\
			\item[d.]  {$G_{m+1}= R(G_{m}^{\prime})$}
		\end{enumerate}
	\end{algorithm}
\end{center}

Let $G_{m}$ be the output of the algorithm \ref{alg13}. Thus if this algorithm terminates on the set $m^{th}$ iteration. The set $G_{m}$ is reduced Gr\"{o}bner basis.
\paragraph{Proof}
\begin{enumerate}[noitemsep]
	\item We first show by induction on $m$ that for each $m^{th}$ iteration every  S-polynomial has a standard representation $S(g_i,g_j,p,q)= \sum_{k=1}^{n} w_k f_k z_k + h_{ij}$. Consider $m=1; G_1=R(\{ f_1, f_2, \dots , f_n \})$; the algorithm produces $f=S(g_i,g_j,p,q)= \sum_{k=1}^{n} w_k f_k z_k + h_{ij}$ as a reduced of $S(g_i,g_j,p,q)$ with respect to $S\cup R$. If $h_{ij} \neq 0$, then $h_{ij} \in G_2$ and again $f$ has a standard representation with respect to $G_2$. Suppose that the hypothesis hold true for $m$. We now prove for $m+1$. If the algorithm terminates at $m+1$ then $G_{m+2}= G_{m+1}=G_{m}$ and hence $f=S(g_i,g_j,p,q)$ has a standard representation with respect to $G_{m+2}= G_{m+1}$. If the algorithm does not terminate at $m+1$, $G_{m+1}=R(G_m \cup \lbrace S(g_i,g_j,p,q) \rbrace)$ so that the algorithm reduces $f=S(g_i,g_j,p,q)$ to $h_{ij}$. This ensures that $f$ has a standard representation with respect to $G_{m+2}$.  Hence the hypothesis hold true for $m+1$. By induction the statement hold true for all $m$.
	\item We now show that the algorithm terminates at $m+1$ if and only if $G_m$ is a finite Gr\"{o}bner basis of $I$: If the algorithm terminates at some $m+1$ then all $S(g_i,g_j,p,q)=0$ and $G_{m+1}=R(G_m)=G_m$, for $G_m$ is a reduced set at every step. Since $\langle G_m \rangle = I$ then we conclude that $G_m$ is a finite reduced Gr\"{o}bner basis for $I$. Conversely if $G_m$ is a finite reduced Gr\"{o}bner basis of $I$, then $R(G_m)=G_m$ and for each pair $(g_i,g_j) \in G_m$, $f=S(g_i,g_j,p,q)$ is reduced to zero by $G_m$. Therefore the  algorithm terminates at $G_{m+1}$. 
	\item If the algorithm never terminates, Let $G=\cup_{m=1}^{\infty} G_m$, then for $m$ sufficiently large every S-polynomial $S(g_i,g_j,p,q)$ has a standard representation with respect to $G_{m+1} \subset G$. Obviously $\langle G \rangle = I$ and hence $G$ is an infinite Gr\"{o}bner basis of $I$.
\end{enumerate}

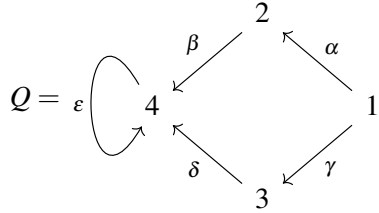
\begin{example}
		\[Q=
		\begin{tikzcd}
		& 2 \arrow[dl,swap,"\beta"] \\
		4 \arrow[out=120,in=240, loop,swap, "\epsilon"] &  & 1 \arrow[ul,swap,"\alpha"] \arrow[dl,"\gamma"] \\
		& 3 \arrow[ul,"\delta"]
		\end{tikzcd}
		\]
		with the length lexicographic ordering $v_1 \prec \dots \prec v_4 \prec \epsilon \prec \beta \prec \delta \prec \alpha \prec \gamma $. Let $f= \alpha \beta - \gamma \delta$, $g= \beta \epsilon$ and $h= \epsilon^{3}$. We see that $LM(f)= \gamma \delta,\quad  LM(g) = \beta \epsilon$ and $LM(h)=  \epsilon^{3}$. $LM(f)\not | LM(h)$ and $LM(f) \not | LM(h)$. The  only possible S-polynomial is $S(g,h,\epsilon^{2}, \beta)= 0$. Thus the set G=$\lbrace f, g ,h \rbrace$ is the Gr\"{o}bner basis since  all the S-polynomial reduces to $0$.  On the other hand if we consider another admissible order $v_1 \prec \dots \prec v_4 \prec \epsilon \prec \beta \prec \delta \prec \gamma \prec \alpha $. We now see that $LM(f)= \alpha \beta,\quad  LM(g) = \beta \epsilon$ and $LM(h)=  \epsilon^{3}$. In this case the only S-polynomial possible is $S(f, g, \epsilon, \alpha)= (\alpha \beta - \gamma \delta)\epsilon - \alpha(\beta \epsilon)=-\gamma \delta \epsilon.$
		$LM(S(f,g,\epsilon,\alpha)) = \gamma \delta \epsilon \notin \langle LM(F), LM(g), LM(h)\rangle$, Thus $G= \lbrace f, g, h \rbrace$ is not a Gr\"{o}bner basis for $I=\langle G \rangle$. We add $r= \gamma \delta \epsilon$ to $G$, and we set $G= \lbrace f, g, h, r \rbrace$. Therefore $S(f,g,\epsilon,\alpha) = r$ and there are no further possible S-polynomial relations. Thus $R(G)=G= \lbrace f, g, h, r \rbrace$ is a Gr\"{o}bner basis for $I$.
\end{example}

In the work of \cite{micah}, the author characterizes quivers whose path algebra has finite Gr\"{o}bner basis. It follows that we can use the above procedures to exhaustively study finite Gr\"{o}bner path algebras.

\end{document}